\documentclass{amsart}
\usepackage{amssymb, amsmath, color}

\usepackage{color}

\newtheorem{theorem}{Theorem}

\newtheorem{corollary}[theorem]{Corollary}
\newtheorem{remark}[theorem]{Remark}

\begin{document}

\title{Ricci Solitons on Lorentzian Manifolds with Large Isometry Groups}
\author{W. Batat $\,$ M. Brozos-V\'{a}zquez $\,$ E. Garc\'{\i}a-R\'{\i}o $\,$ S. Gavino-Fern\'{a}ndez}
\address{WB: \'{E}cole Normale Superieure de L'Enseignement { Technologique} d'Oran \\
{ D\'{e}par\-tement de Math\'{e}matiques et Informatique, B.P. 1523 \\ El M'Naouar, Oran, Algeria}}
\email{wafa.batat@enset-oran.dz}
\address{MBV: Department of Mathematics, University of A Coru\~na, Spain}
\email{mbrozos@udc.es}
\address{EGR-SGF: Faculty of Mathematics,
University of Santiago de Compostela,
15782 Santiago de Compostela, Spain}
\email{eduardo.garcia.rio@usc.es $\,\,$ sandra.gavino@usc.es}
\thanks{Supported by ENSET-Oran (Algeria) and projects MTM2009-07756 and INCITE09 207 151 PR (Spain)}

\begin{abstract}
We show that Lorentzian manifolds whose isometry group is of dimension at least
$\frac{1}{2}n(n-1)+1$ are expanding, steady and shrinking Ricci solitons and
steady gradient Ricci solitons. This provides  examples of complete locally
conformally flat and symmetric Lorentzian Ricci solitons which are not rigid.
\end{abstract}

\maketitle

\section{Introduction}\label{se:1}

A \emph{Ricci soliton} is a pseudo-Riemannian manifold $(M,g)$ which admits
a smooth vector field $X$ on $M$ such that
\begin{equation}\label{soliton}
{\mathcal L}_{X}g+Ric=\lambda g,
\end{equation}
where ${\mathcal L}_{X}$ denotes the Lie derivative in the direction of $X$,
$Ric$ is the Ricci tensor and $\lambda$ is a real number
($\lambda=\frac{1}{n}(2\text{div} X+Sc)$, where $n=\dim\,M$ and $Sc$ denotes the scalar curvature of $(M,g)$).
A Ricci
solition is said to be \emph{shrinking, steady} or \emph{expanding},
 if $\lambda >0,$ $\lambda =0$ or $\lambda <0$, respectively. Moreover we say
that a Ricci soliton $(M,g)$ is a \emph{gradient Ricci soliton} if the
vector field $X$ satisfies $X=\text{grad}\, h$, for some potential function $h$. In such a case equation \eqref{soliton} can be written in terms of $h$ as
\begin{equation}\label{gradsoliton}
2\text{Hes}_{h}+Ric=\lambda g.
\end{equation}
A gradient Ricci soliton is \emph{rigid} if it is isometric to a quotient of $N\times\mathbb{R}^k$,
where $N$ is an Einstein manifold and the potential function is a generalization of the Gaussian soliton (i.e.,  $h=\frac{\lambda}{2}\| x\|^2$ on the Euclidean factor). Rigid solitons have been systematically  studied \cite{PW1} and further investigated in \cite{FL-GR}, \cite{MS}.

Although Ricci solitons exist on many Lie groups and homogeneous spaces, all homogeneous gradient Ricci solitons are rigid \cite{PW2} in the Riemannian setting.
This is based on the existence of splitting results originated by Killing vector
fields on gradient Ricci solitons. Indeed Petersen and Wylie showed that for any Killing vector field $Z$ on a gradient Ricci soliton with potential function $h$, one has that either $Z(h)=0$ or the metric splits off a Euclidean factor. This splitting is not guaranteed in the
pseudo-Riemannian setting due to the fact that $\text{grad} \nabla_Z h$ may be a null vector, in which case $(M,g)$
becomes a Walker manifold (see \cite{walker-metrics} for information on Walker geometry).

For any gradient Ricci soliton with potential function $h$ the following equation is satisfied (see, for instance, \cite{hamilton}):
\[
Sc+2\|\text{grad}\, h\|^2-2\lambda h=C\,,
\]
where $C$ is a constant. Hence, if the scalar curvature is constant, then after rescaling the potential function one may assume
$\|\text{grad}\, h\|^2=\lambda h$, which shows that $\text{grad}\, h$  is a null vector field on $(M,g)$ if and only if the soliton is steady. Further in this
case $Ric(\text{grad}\, h)=0$ (since $Ric(\text{grad}\, h,\,\cdot\,)=dSc(\,\cdot\,)$ holds true for any gradient Ricci soliton) and then one has
that $g(\nabla_{\text{grad}\, h}\text{grad}\, h,Y)=0$ for all vector fields $Y$, thus showing that $\text{grad}\, h$ is a geodesic vector field.
On the other hand note that the gradient of the potential function is a recurrent vector field (i.e., the plane field $\text{span}\{\text{grad}\, h\}$ is parallel)
if and only if $\nabla_Y\text{grad}\, h=\omega(Y)\text{grad}\, h$ for some $1$-form $\omega$. This means that the $(1,1)$-Hessian tensor
$hes_h(Y)=\nabla_Y\text{grad}\, h=\omega(Y)\text{grad}\, h$ for all vector fields $Y$ and thus the constancy of the scalar curvature gives
$Ric(\text{grad}\, h)=-hes_h(\text{grad}\, h)=0$, showing that if the gradient of the potential function is recurrent then the Ricci operator and the Hessian tensor are two-step nilpotent.

The situation described above occurs, for instance, when studying pseudo-Rie\-man\-nian manifolds with large isometry groups.
Riemannian manifolds admitting a group of isometries of dimension at least
$\frac{1}{2}n(n-1)+1$ are either of constant curvature or products
of an $(n-1)$-dimensional space of constant curvature with a line or a circle.
Lorentzian metrics allow other solutions which are related to the existence
of null submanifolds on $M$ which are left invariant by Lie group actions (see, for instance, \cite{BV-C-GR-GF}).
If the group of isometries has dimension at least $\frac{1}{2}n(n-1)+2$, then the sectional
curvature is constant \cite{P}. However, if the group of isometries has dimension
$\frac{1}{2}n(n-1)+1$, then non-homogeneous examples exist.
The complete classification of Lorentzian manifolds
with an isometry group of this large was given for any dimension $n \geq 4$, $n\neq 7$, by Patrangenaru \cite{P}. Besides the spaces of constant curvature
$M_{1}^{n}(c)$ and manifolds reducible as products
$N^{n-1}(c)\times \mathbb{R}$, the remaining examples are:
\begin{itemize}
\item \emph{Egorov spaces:} Lorentzian manifolds $(\mathbb{R}^{n+2},g_{f})$,
where $f$ is a positive function of a real variable and
\begin{equation}\label{ego}
g_{f}(u,v,x_1,\dots,x_n)=du\,dv+f(u)\sum_{i=1}^{n}(dx_{i})^{2}.
\end{equation}
\end{itemize}
\begin{itemize}
\item \emph{$\varepsilon$-spaces:} Lorentzian manifolds $(\mathbb{R}^{n+2},g_{\varepsilon })$, where
\begin{equation}\label{epsi}
g_{\varepsilon }(u,v,x_1,\dots,x_n)=\varepsilon
\sum_{i=1}^{n}x_{i}^{2}(du)^{2}+du\,dv+\sum_{i=1}^{n}(dx_{i})^{2}.
\end{equation}
\end{itemize}

The geometry of Egorov spaces and $\varepsilon$-spaces has been investigated in
the literature (see for example
\cite{BCD}, \cite{clptv}, \cite{CG}, \cite{E}, \cite{P2} and references therein)
where it is shown that both are Walker manifolds.
It is worth emphasizing here that although $\varepsilon$-spaces are locally symmetric,
Egorov spaces are not homogeneous in general. However both families are locally conformally flat
and have two-step nilpotent Ricci operator.

Recall that spaces of constant curvature are Einstein and thus trivial Ricci solitons
and moreover, products $N^{n-1}(c)\times\mathbb{R}$ are rigid Ricci solitons.
The first purpose of this note is to show that both Egorov spaces and $\varepsilon$-spaces are (expanding, steady and shrinking) Ricci solitons.
Secondly we show that both families are gradient Ricci solitons, that are necessarily steady, and therefore non-rigid.
This leads to new examples of Lorentzian locally conformally flat Ricci solitons without Riemannian analog
(see \cite{PW3} for the Riemannian case).

This work is structured as follows. Egorov spaces are studied in Section \ref{se:2}. In Section \ref{se:3} we consider the Cahen-Wallach Lorentzian symmetric spaces that describe the indecomposable but not irreducible
Lorentz symmetric spaces \cite{cahen-wallach}, \cite{clptv}, and that generalize $\varepsilon$-spaces. We show that they all are expanding, steady and shrinking
Ricci solitons but only steady Ricci solitons may be gradient ones. The $\varepsilon$-spaces are then obtained as the locally conformally flat Cahen-Wallach Lorentzian symmetric spaces.

\section{Egorov spaces}\label{se:2}
Let $(\mathbb{R}^{n+2},g_{f})$, $n\geq 1$ denote an Egorov space. As proved
in \cite{BCD}, with respect to the basis of coordinate vector fields
$\{\partial _{u}=\frac{\partial }{\partial {u}}$,
$\partial _{v}=\frac{\partial}{\partial {v}}$,
$\partial _{i}=\frac{\partial }{\partial {x_{i}}}\}$, with $i=1,\dots ,n$,
for which $g_f$ adopts expression \eqref{ego}, the non-vanishing covariant derivatives of coordinate
vector fields are given by
\begin{equation}\label{nablaei}
\nabla _{\partial _{i}}\partial _{i}=-\frac{f^{\prime }}{2}\partial_{v},
\qquad
\nabla _{\partial _{i}}\partial _{u}=\frac{f^{\prime }}{2f}\partial _{i},
\qquad i=1,\dots ,n.
\end{equation}
Hence observe that $\partial_v$ is a parallel null vector field and thus that Egorov spaces
are Walker metrics \cite{CG}, \cite{walker-metrics}.
By an explicit calculation on the geodesic equations, it has been shown in \cite{BCD} that
Egorov spaces are geodesically complete.
The curvature tensor $R$, which is given by
$R(X,Y)=\nabla _{[X,Y]}-[\nabla_{X},\nabla_{Y}]$,
is determined by
\begin{equation}\label{curvei}
R_{iui}^{v}=\frac{1}{4f}\left[ (f^{\prime })^{2}-2ff^{\prime \prime }\right],
\quad
R_{iuu}^{i}=-\frac{1}{4f^{2}}\left[ (f^{\prime })^{2}-2ff^{\prime
\prime }\right],
\quad i=1,\dots ,n.
\end{equation}
The Ricci tensor
$Ric\left( X,Y\right)=\mathrm{trace}\left\{ Z\rightarrow R\left( X,Z\right) Y\right\}$
satisfies
\begin{equation}\label{rhoei}
Ric_{uu}=\frac{n}{4f^{2}}[(f^{\prime })^{2}-2ff^{\prime \prime }],
\end{equation}
being zero otherwise. This shows that the Ricci operator is two-step nilpotent.

\begin{remark}\rm
In opposition to $\varepsilon$-spaces, Egorov spaces are not homogeneous in general.
However the Ricci tensor is recurrent and so is the curvature tensor since they are locally
conformally flat (see \cite{BCD}, \cite{CG}).
\end{remark}

\subsection{Soliton equations}\label{se:2.1}

Next we will show that all Egorov metrics are Ricci solitons.
Let $\displaystyle X=\sum_{l=u,v,1,..,n}X_{l}\,\partial_{l}$
be an arbitrary vector field on $(\mathbb{R}^{n+2},g_{f})$. Then
\eqref{soliton} becomes
\begin{equation}\label{syssoliego}
\left\{
\begin{array}{lc}
\partial _{i}X_{j}+\partial _{j}X_{i}=0,&\quad 1\leq i\neq j\leq n, \\
\noalign{\medskip}
\partial _{i}X_{v}+f\partial _{u}X_{i}=0,&\quad 1\leq i\leq n, \\
\noalign{\medskip}
\partial _{i}X_{u}+f\partial _{v}X_{i}=0,&\quad 1\leq i\leq n, \\
\noalign{\medskip}
X_{u}f^{\prime }+2f\partial _{i}X_{i}=\lambda f,&\quad 1\leq i\leq n, \\
\noalign{\medskip}
Ric _{uu}+2\partial _{u}X_{v}=0,&\quad\\
\noalign{\medskip}
\partial _{u}X_{u}+\partial _{v}X_{v}=\lambda ,
\quad
\partial _{v}X_{u}=0. &\quad
\end{array}%
\right.
\end{equation}
Now it follows by a straightforward calculation that the metric $g_{f}$ is a Ricci soliton since the vector
field
\begin{equation}\label{eq:part-solution-egorov}
X=\left( -\frac{1}{2}\int Ric _{uu}\mathrm{d}u+\lambda v\right) \partial _{v}
+
\underset{i=1}{\overset{n}{\sum }}\frac{\lambda }{2}\,x_{i}\, \partial _{i},
\end{equation}
satisfies \eqref{syssoliego}. Note that $\lambda$ is the constant of equation \eqref{soliton} and can be chosen with absolute freedom. Therefore, we obtain the following:

\begin{theorem}\label{th:2}
All Egorov spaces $(\mathbb{R}^{n+2},g_{f})$ are expanding, steady and shrinking Ricci solitons.
\end{theorem}

\begin{remark}\label{remark:general-solution-egorov}\rm
By a standard process of integration, one gets that the general solution of  \eqref{syssoliego}, is
\[
\begin{array}{rcl}
X_u&=& a+b\, u\,, \\
\noalign{\bigskip}
X_v&=&c_0+ (\lambda-b) v-\frac{1}{2}\int Ric_{uu}\mathrm{d}u+\sum_i k_i x_i
+(\sum_i x_i^2)K\,,\\
\noalign{\bigskip}
X_i&=&c_i-\int \frac{k_i}{f}\,\mathrm{d}u+\left(\frac{\lambda}{2}-(a+bu)\frac{f'(u)}{2 f(u)}\right) x_i+\sum_{j\neq i} A_{ij} x_j\,,
\end{array}
\]
where $a$, $b$ and $K$ are constants which satisfy the equation
\[
b f'(u)+(a+b u)\left(f''(u)-\frac{f'(u)^2}{f(u)}\right)=4K,
\]
$(A_{ij})$ is an arbitrary skew-symmetric matrix, $c_0$, $c_i$ and $k_i$ are arbitrary constants for $i=1,\dots,n$.
\end{remark}

\begin{remark}\rm
The general solution of \eqref{syssoliego} obtained in Remark~\ref{remark:general-solution-egorov} can be determined from the particular solution \eqref{eq:part-solution-egorov} using the following observation. Any two vector fields $X_1$ and $X_2$ satisfying \eqref{soliton}
($\mathcal{L}_{X_i}g+Ric=\lambda_ig$, $i=1,2$) differ in a conformal vector field with constant divergence (i.e., a homothetic vector field) since
\[
\mathcal{L}_{X_1-X_2}g-(\lambda_1-\lambda_2)g =\mathcal{L}_{X_1}g-\lambda_1g -\mathcal{L}_{X_2}g+\lambda_2g
=0.
\]
Conversely, for any Ricci soliton $X_1$, adding a homothetic vector field gives another Ricci soliton. A special case of the above occurs if the manifold is compact, where an immediate application of the divergence theorem shows that any two Ricci solitons differ in a Killing vector field.
\end{remark}

\subsection{Gradient Ricci solitons}\label{se:2.2}

Now, let $X=\mathrm{grad}\,h$ be an arbitrary gradient vector field on $(\mathbb{R}^{n+2},g_{f})$ with potential function $h$ (which is given by
$\text{grad}\, h=\left(\partial_v h,\partial_u h,\frac{1}{f}\partial_{1}h,\dots,\frac{1}{f}\partial_{n}h  \right)$ in the coordinates \eqref{ego}).
By standard calculations we get from \eqref{syssoliego}  that, $(\mathbb{R}^{n+2},g_{f})$ is a gradient Ricci soliton if and only if the following holds
\begin{equation}\label{sysgradego}
\left\{
\begin{array}{lc}
f^{\prime }\partial _{v}h+2\partial^2_{ii}h=\lambda f, &\quad 1\leq i\leq n, \\
\noalign{\medskip}
2\partial^2_{iu}h-\frac{f^{\prime }}{f}\partial_{i}h=0, &\quad 1\leq i\leq n, \\
\noalign{\medskip}
2\partial^2_{uu}h+Ric _{uu}=0, &\quad  \\
\noalign{\medskip}
\partial^2_{uv}h=\frac{\lambda }{2}, &\quad \\
\noalign{\medskip}
\partial^2_{ij}h=\partial^2_{iv}h=\partial^2_{vv}h=0, &\quad 1\leq i\neq j\leq n.
\end{array}
\right.
\end{equation}
Integrating the equations $\partial^2_{iv}h=\partial^2_{vv}h=0$ in \eqref{sysgradego} we obtain that the potential function
splits as $h(u,v,x_1,\dots,x_n)=h_0(u,x_1,\dots,x_n)+h_1(u)v$ for some functions $h_0$, $h_1$.
Moreover equations $\partial^2_{uv}h=\frac{\lambda }{2}$ and
 $\partial^2_{ij}h=0$ now show that
$h(u,v,x_1,\dots,x_n)=h_0(u)+\sum_i h_i(u,x_i) +(\frac{\lambda}{2}u+\kappa)v$ for some constant $\kappa$
and functions $h_i$, $i=1,\dots,n$. Hence \eqref{sysgradego} reduces to
\begin{equation}\label{sysgradego2}
\left\{
\begin{array}{lc}
(\kappa+\frac{\lambda}{2}u)f^\prime + 2 \partial^2_{ii}h_i=\lambda f, &\quad 1\leq i\leq n, \\
\noalign{\bigskip}
2 f\partial^2_{iu}h_i=f^\prime \partial_{i}h_i, &\quad 1\leq i\leq n, \\
\noalign{\medskip}
Ric_{uu}+2h_0^{\prime\prime}+2\sum_i\partial^2_{uu}h_i=0. &\quad
\end{array}
\right.
\end{equation}
Integrating the first equations in \eqref{sysgradego2} one gets
\[
h_i(u,x_i)=\frac{1}{8} x_i^2\left( 2\lambda f-(2\kappa+u\lambda)f^\prime \right)+x_i h_i(u)+k_i(u)
\]
for some functions $h_i(u)$ and $k_i(u)$. Substituting the above into \eqref{sysgradego2}
and deriving the second equation
$2f\partial^3_{iiu}h_i=f^\prime \partial^2_{ii}h_i$ one gets
\[
(\frac{\lambda}{2}u+\kappa)\left((f^{\prime })^{2}-2ff^{\prime\prime}\right)=0\,,
\]
which shows that $\lambda=\kappa=0$ unless the manifold is flat. Hence the potential function becomes
$h(u,v,x_1,\dots,x_n)=h_0(u)+\sum_i x_i h_i(u)$. Now the second equation in \eqref{sysgradego2}
reduces to
\[
2 f h_i^\prime = f^\prime h_i
\]
and hence $h_i(u)=c_i\sqrt{f(u)}$ for some constants $c_i$.
Now, deriving the last equation in \eqref{sysgradego2}
with respect to $x_i$, one gets
$0=c_i ((f^\prime)^2-2ff^{\prime\prime})$ and thus all $c_i$ must vanish. Therefore, the potential
function reduces to a function of the single variable $u$,
$h(u,v,x_1,\dots,x_n)=h_0(u)$, and hence the only remaining constraint in  \eqref{sysgradego2}
is $Ric_{uu}+2h_0^{\prime\prime}=0$.
Therefore, we have shown that

\begin{corollary}
A non-flat Egorov space is a gradient steady Ricci soliton and the potential function $h(u)$ is given by $h^{\prime\prime}=-\frac{1}{2}Ric_{uu}$.
\end{corollary}

\begin{remark}\rm\label{re:gsrs1}
Note that for any function $h(u)$, $\text{grad}\, h=h^\prime(u)\partial_v$ is a null vector field and moreover $\nabla_{\partial_u}\text{grad}\, h=h^{\prime\prime}(u)\partial_v$ (the other derivatives being zero).
Therefore, $\text{grad}\, h$ is a null geodesic vector field.
Further observe that $\text{grad}\, h$ is a recurrent vector field in the direction of the parallel null vector $\partial_v$.
\end{remark}


\section{Cahen-Wallach symmetric spaces}\label{se:3}

Recall that the notion of irreducibility (the holonomy group does not stabilize any nontrivial subspace)
is very strong in the pseudo-Riemannian setting. Indeed, irreducible Lorentzian symmetric spaces
are necessarily of constant sectional curvature \cite{clptv}. A pseudo-Riemannian manifold is
said to be \emph{indecomposable} if the holonomy group, acting at each point $p\in M$, stabilizes only non-trivial degenerate subspaces $V\subset T_pM$ (i.e., the restriction of the metric to $V\times V$
is degenerate).

Clearly any irreducible Lorentzian symmetric space is a trivial (Einstein) Ricci soliton. Our purpose in
this section is to show that indecomposable (not irreducible) Lorentzian symmetric spaces are non-trivial Ricci solitons.
Indecomposable Lorentzian symmetric spaces are either irreducible or the so-called
Cahen-Wallach symmetric spaces which are given as follows \cite{cahen-wallach}, \cite{clptv}.
Take  $M=\mathbb{R}^{n+2}$  and define a metric tensor by
\begin{equation}\label{CW-metric}
g_{cw}(u,v,x_1,\dots,x_n)=\left(\sum_{i=1}^n \kappa_i x_i^2\right) (du)^2+du\, dv
+\sum_{i=1}^n (dx_i)^2\,,
\end{equation}
where $\kappa_i$, $i=1,\dots,n$, are non-zero real numbers.
The Levi-Civita connection is determined by the non-zero  Christoffel symbols:
\[
\nabla_{\partial_u} \partial_u = -\sum_i \kappa_i x_i \partial_{i},
\qquad
\nabla_{\partial_u} \partial_{i} = \nabla_{\partial_{i}} \partial_u= \kappa_i x_i \partial_v\,.
\]
Hence, the only non-zero components of the $(0,4)$ curvature tensor
are given by
\[
R(\partial_u,\partial_{i},\partial_u,\partial_{i})=-\kappa_i,
\quad i=1,\dots,n.
\]
The Ricci tensor satisfies
$Ric_{uu}=-\sum_{i=1}^n\kappa_i$, the other terms being zero.
Hence the Ricci operator is $2$-step nilpotent (and thus the scalar curvature vanishes).

Cahen-Wallach symmetric spaces are locally conformally flat if and only if $\kappa_1=\cdots=\kappa_n$,
in which case the resulting manifold is an $\varepsilon$-space (we refer to \cite{cahen-wallach} and \cite{clptv} for more information
on the geometry of Cahen-Wallach spaces and to \cite{BCD} and \cite{CG} for $\varepsilon$-spaces).

\subsection{Soliton equations}\label{se:3.1}

Let $X=X_u \partial_u+X_v \partial_v+\sum_i X_i \partial_{i}$ be an arbitrary vector field on
$\mathbb{R}^{n+2}$. The Lie derivative of the metric $g$ is given by
\[
\mathcal{L}_Xg\!=\!\left(\!\!\begin{array}{ccc}
2(\sum_i \kappa_i x_i X_i+(\sum_i \kappa_i x_i^2)\partial_u X_u+\partial_u X_v)
& A_{uv} & A_{uj}\\
\noalign{\medskip}
A_{uv} & 2\partial_v X_u & \partial_{j}X_u+\partial_v X_j \\
\noalign{\medskip}
A_{uj} & \partial_{j}X_u+\partial_v X_j & A_{jk}
\end{array}\!\!\right),
\]
where
\[
\begin{array}{l}
A_{uv}=\partial_u X_u+\partial_v X_v+(\sum_i \kappa_i x_i^2)\partial_v X_u,\\
\noalign{\medskip}
A_{uj}=(\sum_i \kappa_i x_i^2)\partial_{j} X_u+\partial_{j} X_v+\partial_u X_j, \,\, \text{and}\\
\noalign{\medskip}
A_{jk}=\partial_{j}{X_k}+\partial_{k}{X_j}.
\end{array}
\]
Thus the Ricci soliton equation \eqref{soliton} is equivalent to (note that $\lambda=\frac{2}{n+2} \text{div}(X)$
and $\text{div}(X)=\partial_u{X_u}+\partial_v{X_v}+\sum_i \partial_{i}{X_i}$):
\begin{equation}\label{eq:equation-soliton}
\left\{ \begin{array}{lcl}
\sum_i \kappa_i-2\sum_i \kappa_i x_i X_i-2\partial_u X_v
-(\sum_i \kappa_i x_i^2) (2 \partial_u{X_u}
-\lambda)
&=&0\\
\noalign{\bigskip}
\partial_v{X_u}&=&0\\
\noalign{\bigskip}
(\sum_i \kappa_i x_i^2) \partial_v{X_u}+(\partial_u{X_u}+\partial_v{X_v})
-\lambda&=&0\\
\noalign{\bigskip}
(\sum_i \kappa_i x_i^2)\partial_{j} X_u
+\partial_{j} X_v+\partial_u X_j&=&0\\
\noalign{\bigskip}
\partial_{j}{X_u}+\partial_v{X_j}&=&0\\
\noalign{\bigskip}
2\partial_{j}{X_j} -\lambda&=&0\\
\noalign{\bigskip}
\partial_{j}{X_k}+\partial_{k}{X_j}&=&0
\end{array}\right.
\end{equation}
where  $j,k=1,\dots,n$ and $j\neq k$.

First of all, observe that a particular solution of the previous system is:
\begin{equation}\label{particular-1}
X_{(u,v,x_1,\dots,x_n)}=\left( 0, \frac{1}{2}(\sum_i \kappa_i)u+\lambda v,\frac{\lambda}{2} x_1,\dots,\frac{\lambda}{2} x_n\right) ,
\end{equation}
where $\lambda$ is the constant of equation \eqref{soliton} and can be chosen with absolute freedom. The choice of this particular family of spacelike vector fields shows that $(\mathbb{R}^{n+2},g_{cw})$ is a expanding, steady
or shrinking Ricci soliton, and therefore proves the following:

\begin{theorem}\label{th:4}
Indecomposable Lorentzian symmetric spaces are expanding, steady and shrinking Ricci solitons.
\end{theorem}

\begin{remark}\label{remark:general-solution-cw}\rm
Next we will consider the general solution of \eqref{eq:equation-soliton} assuming that
$\kappa_i\neq\kappa_j$ for some $i,j\in\{ 1,\dots,n\}$
(the case $\kappa_1=\cdots=\kappa_n$ will be consider in Section \ref{se:3.2}).
After a standard integration process one gets that all Ricci solitons are given by
vector fields $X=(X_u,X_v,X_1,\dots,X_n)$ of the form
\[
\begin{array}{l}
X_{u}=a,\\
\noalign{\bigskip}
X_v=b+\lambda v+(\sum_i \kappa_i)\frac{u}{2}-\sum_i x_i h_i'(u),
\\
\noalign{\bigskip}
X_{j}=\frac{\lambda}{2}x_j+h_j(u)+\sum_{i\neq j} c_{ij} x_i,\quad j=1,\dots, n,
\end{array}
\]
where $a$, $b$, are arbitrary constants, $c_{ij}$ is a skew-symmetric matrix of constants satisfying $c_{ij}(\kappa_i-\kappa_j)=0$ and $h_i$ are functions which satisfy $\kappa_i h_i-h_i''=0$ (hence note that if $\kappa_i>0$ then $h_i(u)=d_i^1 e^{-u\sqrt{\kappa_i}} +d_i^2 e^{u\sqrt{\kappa_i}}$ while, if $\kappa_i<0$ then  $h_i(u)=d_i^1 \sin(u\sqrt{-\kappa_i}) +d_i^2 \cos(u\sqrt{-\kappa_i})$ for  $d_i^{1,2}\in \mathbb{R}$).

A straightforward calculation shows that the causal character of $X$ depends on the point and the value of $a$, $b$, $c_{ij}$, $h_j$, $\kappa_i$ and $\lambda$ ($i,j=1,\dots,n$).
\end{remark}

\subsection{Gradient Ricci solitons}\label{se:3.2}

Let $h$ be a function on $\mathbb{R}^{n+2}$. Then the gradient with respect to the metric
\eqref{CW-metric} is given by
\[
\text{grad}(h)=(\partial_vh,-(\sum_i \kappa_i x_i^2) \partial_vh+\partial_uh,\partial_1h,\dots,\partial_jh,\dots,\partial_nh)\,.
\]
and thus (\ref{eq:equation-soliton}) becomes (where  $j,k=1,\dots,n$, $j\neq k$)
\begin{equation}\label{eq:equation-gradient-soliton}
\left\{ \begin{array}{lcl}
\sum_i \kappa_i-2\sum_i \kappa_i x_i \partial_{i}{h}-2\partial^2_{uu} h
+(\sum_i \kappa_i x_i^2) \lambda&=&0\\
\noalign{\bigskip}
\partial^2_{vv} h&=&0\\
\noalign{\bigskip}
2\partial_{uv}^2h-\lambda&=&0\\
\noalign{\bigskip}
\kappa_j x_j \partial_v{h}-\partial^2_{uj} h&=&0\\
\noalign{\bigskip}
\partial^2_{vj} h&=&0\\
\noalign{\bigskip}
2\partial^2_{jj} h-\lambda&=&0\\
\noalign{\bigskip}
\partial^2_{jk} h&=&0
\end{array}\right.
\end{equation}
Using that
$\partial^2_{vv} h$ $=$ $\partial^2_{vj} h$ $=$ $\partial^2_{jk} h=0$, $(j\neq k)$, the function $h$ separates variables $v,x_1,\dots,x_n$ and
the previous system  reduces to
\begin{equation}\label{eq:sist-interm}
\left\{ \begin{array}{lcl}
\sum_i \kappa_i-2\sum_i \kappa_i x_i \partial_{i}{h}-2\partial^2_{uu} h
+(\sum_i \kappa_i x_i^2) \lambda&=&0\\
\noalign{\bigskip}
2\partial_{uv}^2h-\lambda&=&0\\
\noalign{\bigskip}
\kappa_j x_j \partial_v{h}-\partial^2_{uj} h&=&0\\
\noalign{\bigskip}
2\partial^2_{jj} h-\lambda&=&0
\end{array}\right.
\end{equation}
where  $j=1,\dots,n$.

Differentiate the third equation with respect to $u$ to obtain that $\kappa_j x_j \partial^2_{uv} h-\partial^3_{uuj} h=0$, and, by the second equation we have
\begin{equation}\label{eq:paso-interm}
2\partial^3_{uuj} h=\kappa_j x_j \lambda\,.
\end{equation}
Now, differentiate the first equation with respect to $x_j$ to obtain
\[
-2\kappa_j\partial_{j} h-2\kappa_j x_j \partial^2_{jj} h+2\kappa_j x_j \lambda -2 \partial^3_{uuj} h=0\,.
\]
Use equations in \eqref{eq:sist-interm} and \eqref{eq:paso-interm} to simplify and get that
\[
\kappa_j\partial_{j} h=0\,.
\]
This shows that $h$ does not depend on the variables $x_j$,  $(j=1,\dots, n)$. Simplify now the third equation of (\ref{eq:sist-interm}) to conclude that $\partial_v h=0$, so $h$ is only a function of the variable $u$, $h=h(u)$, and (\ref{eq:equation-gradient-soliton})  reduces to
\[
\sum_i \kappa_i-2h^{\prime\prime}(u)=0\,.
\]
We integrate this equation to obtain the following function as the general solution of the system of equations (\ref{eq:equation-gradient-soliton}):
\[
h(u)=\alpha + \beta u + \frac{1}{4}\sum_i \kappa_i u^2\,.
\]
Note that $\text{grad}(h)$ is a null vector field that satisfies equation \eqref{soliton} for $\lambda=0$, indeed $\lambda=\frac{2div(grad(h))}{n+2}=\Delta h=0$. Therefore, we have proved the following

\begin{theorem}\label{th:5}
Indecomposable Lorentzian symmetric space are steady gradient Ricci solitons.
\end{theorem}

\begin{remark}\rm
Finally observe that, proceeding as in Remark \ref{re:gsrs1},
$X=grad(h)=\left(\beta+\frac{1}{2}u\sum_i \kappa_i\right)\partial_v$ is a null geodesic vector field which is recurrent. Hence, the line field $L=\text{span}\{\text{grad}(h)\}$ is a parallel degenerate one-dimensional plane field, which agrees with the Walker structure of Cahen-Wallach spaces.
\end{remark}

\begin{remark}\rm
Riemannian complete shrinking Ricci solitons with bounded curvature can be made gradient by adding a Killing vector field \cite{Naber}.
Theorem~\ref{th:5} and Theorem~\ref{th:4} (for a suitable choice of $f$ which makes the curvature bounded) show that this
result does not hold in the Lorentzian setting due to the existence of complete shrinking
or expanding Ricci solitons which cannot be made gradient by adding any Killing vector field.
\end{remark}

\subsection{$\epsilon$-spaces}\label{se:3.3}

Locally conformally flat Cahen-Wallach symmetric spaces are precisely the $\varepsilon$-spaces introduced
in Section \ref{se:1}. In this case the metric is given by \eqref{CW-metric} with
$\kappa_1=\dots=\kappa_n=\varepsilon$.

For $X=(X_u,X_v,X_1,\dots,X_n)$ we particularize the solutions in Remark~\ref{remark:general-solution-cw}
to obtain
\[
\begin{array}{l}
X_{u}=a,\\
\noalign{\bigskip}
X_v=b+\lambda v+n\varepsilon\frac{u}{2}-(\sum_i x_i h'_i(u)),
\\
\noalign{\bigskip}
X_{j}=\frac{\lambda}{2}x_j+h_j(u)+\sum_{i\neq j} c_{ij} x_i,\quad j=1,\dots, n,
\end{array}
\]
where $a$, $b$, are arbitrary constants, $c_{ij}$ is a skew-symmetric matrix of constants  and $h_i$ are functions satisfying $\varepsilon h_i-h_i''=0$.
Therefore \emph{$\varepsilon$-spaces are expanding, steady and shrinking Ricci solitons which
are symmetric and  locally conformally flat}.

Moreover, Theorem \ref{th:5} shows that
\emph{$\varepsilon$-spaces are steady gradient Ricci solitons} with potential function
$h(u)=\alpha + \beta u + \frac{n}{4}\varepsilon u^2$.


\begin{thebibliography}{99}
\bibitem{BCD}
W. Batat, G. Calvaruso, B. De Leo,
Curvature properties of
Lorentzian manifolds with large isometry group,
\emph{Math. Phys. Anal. Geom.} 12 (2009), 201--217.

\bibitem{BV-C-GR-GF}
M. Brozos-V\'{a}zquez, G. Calvaruso, E. Garc\'{\i}a-R\'{\i}o, S. Gavino-Fern\'{a}ndez,
Three-dimensional Lorentzian homogeneous  Ricci solitons, Israel J. Math., to appear.

\bibitem{walker-metrics}
M. Brozos-V\'{a}zquez, E. Garc\'{i}a-R\'{i}o, P. Gilkey, S. Nik\v{c}evi\'{c}, and R. V\'{a}zquez-Lorenzo,
\emph{The geometry of Walker manifolds},
Synthesis Lectures on Mathematics and Statistics \textbf{5}, Morgan \& Claypool Publ., 2009.

\bibitem{cahen-wallach}
M. Cahen, N. Wallach,
Lorentzian symmetric spaces,
\emph{Bull. Amer. Math. Soc.} \textbf{76} (1970), 585--591.

\bibitem{clptv}
M. Cahen, J. Leroy, M. Parker, F. Tricerri, L. Vanhecke,
Lorentz manifolds modelled on a Lorentz symmetric space,
\emph{J. Geom. Phys.} \textbf{7} (1990), 571--581.

\bibitem{CG}
G. Calvaruso, E. Garc\'ia-R\'io,
Algebraic properties of curvature
operators in Lorentzian manifolds with large isometry groups,
\emph{SIGMA Symmetry Integrability Geom. Methods Appl.} \textbf{6}  (2010), Paper 005, 8 pp.

\bibitem{E}
I.P. Egorov,
Riemannian spaces of the first three lacunary types
in the geometric sense,
\emph{Dokl. Akad. Nauk. SSSR} \textbf{150} (1963), 730--732.

\bibitem{FL-GR}
M. Fern\'{a}ndez-L\'{o}pez, E. Garc\'{\i}a-R\'{\i}o,
Rigidity of shrinking Ricci solitons,
\emph{Math. Z.}, to appear.

\bibitem{hamilton}
R. S. Hamilton, The formation of singularities in the Ricci flow, Surveys in Differential Geometry (Cambridge, MA, 1993), Vol. II, 7--136, International Press, Cambridge, MA, 1995.

\bibitem{MS}
O. Munteanu, N. Sesum,
On gradient Ricci solitons,
arXiv:0910.1105v1.

\bibitem{Naber}
A. Naber,
Noncompact shrinking four solitons with nonnegative curvature,
arXiv:0710.5579v2.

\bibitem{onda}
K. Onda,
Lorentz Ricci Solitons on 3-dimensional Lie groups
\emph{Geom. Dedicata}, DOI 10.1007/s10711-009-9456-0, to appear.

\bibitem{P}
V. Patrangenaru,
Lorentz manifolds with the three largest degrees of symmetry,
\emph{Geom. Dedicata} \textbf{102} (2003), 25--33.

\bibitem{P2}
V. Patrangenaru,
Locally homogeneous pseudo-Riemannian manifolds,
\emph{J. Geom. Phys.} \textbf{17} (1995), 59--72.

\bibitem{PW1} P. Petersen and W. Wylie;
{Rigidity of gradient Ricci solitons,}
\emph{Pacific J. Math.} \textbf{241} (2009), 329--345.

\bibitem{PW2} P. Petersen and W. Wylie;
{On gradient Ricci solitons with symmetry},
\emph{Proc. Amer. Math. Soc.}  \textbf{137} (2009), 2085--2092.

\bibitem{PW3} P. Petersen and W. Wylie;
{On the classification of gradient Ricci solitons},
arXiv:0712.1298v5.
\end{thebibliography}
\end{document}